# RIEMANN ZETA FUNCTION

DORIN GHISA*

(*) York University, Glendon College, 2275 Bayview Avenue, Toronto, Canada, M4N 3M6
E-mail address: dghisa@yorku.ca

The domain of the Riemann Zeta function and that of its derivative appear as branched covering Riemann surfaces $(\overline{\mathbb{C}},f)$. The fundamental domains, which are the leafs of those surfaces are revealed. For this purpose, pre-images of the real axis by the two functions are taken and a thorough study of their geometry is performed. The study of intertwined curves generated in this way, allowed us to prove that the Riemann Zeta function has only simple zeros and finally that the Riemann Hypothesis is true. Color visualization of the conformal mappings of the fundamental domains is used for a better understanding of the theory.

## 1. Introduction

The Riemann Zeta function is one of the most studied transcendental functions in view of its many applications in number theory, algebra, complex analysis, statistics, as well as in physics. Another reason why this function has drawn so much attention is the celebrated Riemann Hypothesis (RH) regarding its non trivial zeros, which has resisted proof or disproof until now.

Hopefully, starting from now, people can write RP instead of RH, meaning the Riemann Property (of non trivial zeros).

The RP proof will be derived from the global mapping properties of Zeta function. The Riemann conjecture prompted the study of at lest local mapping properties in the neighborhood of non trivial zeros. There are known color visualizations of the module, the real part and the imaginary part of Zeta function at some of those points (see Wolfram MathWorld), however they do not offer an easy way to visualize the global behavior of this function. We perfected an idea found in [8], page 213.

The Riemann Zeta function has been obtained by analytic continuation [1], page 178 of the series

$$\zeta(s) = \sum_{n=1}^{\infty} n^{-s}, \qquad s = \sigma + it \qquad \#$$

which converges uniformly on the half plane $\sigma \geq \sigma_0$, where $\sigma_0 > 1$ is arbitrarily chosen. It is known [1], page 215, that Riemann function $\zeta(s)$ is a meromorphic function in the complex plane having a single simple pole at $s = 1$ with the residue 1. The representation formula

$$\zeta(s) = -\frac{\Gamma(1-s)}{2\pi i}\int_C [(-z)^{s-1}/(e^z - 1)]dz \qquad \#$$

where $\Gamma$ is the Euler function and $C$ is an infinite curve turning around the origin, which does not enclose any multiple of $2\pi i$, allows one to see that $\zeta(-2m) = 0$ for every positive integer $m$ and that there are no other zeros of $\zeta$ on the real axis. However, the function $\zeta$ has infinitely many other zeros (so called, *non trivial* ones), which are all situated in the *(critical) strip* $\{s = \sigma + it : 0 < \sigma < 1\}$. The famous RH says that these zeros are actually on the *(critical) line* $\sigma = 1/2$.

We will make reference to the Laurent expansion of $\zeta(s)$ for $|s - 1| > 0$:

$$\zeta(s) = \frac{1}{s-1} + \sum_{n=0}^{\infty}[(-1)^n/n!]\gamma_n(s-1)^n, \qquad \#$$

where $\gamma_n$ are the Stieltjes constants:

$$\gamma_n = \lim_{m\to\infty}[\sum_{k=1}^{m}(\log k)^n/k - (\log m)^{n+1}/(m+1)] \qquad \#$$

as well as to the functional equation [1], page 216:

$$\zeta(s) = 2^s \pi^{s-1} \sin\frac{\pi s}{2}\Gamma(1-s)\zeta(1-s). \qquad \#$$

As predicted, the proof of RH will have a lot of consequences in number theory, although RP is a property of conformal mappings. It is probably safe to say that the conformal mappings made their royal entrance into the realm of number theory through the RP. By the global mapping properties of the Riemann Zeta function we understand the conformal mapping properties of its fundamental domains.

The concept of a fundamental domain, as it appears here, has been formulated by Ahlfors, who also emphasized its importance as a tool in the study of analytic functions. He said: *Whatever the advantage of such a representation may be, the clearest picture of the Riemann surface is obtained by direct consideration of the fundamental regions in the z-plane* ([1], page 99).

Starting with the study of Blaschke products ([4], [5], [7]), we followed a true program of revealing fundamental domains for different classes of analytic functions ( [2],[3],[6] ), drawing in the end the conclusion that this can be done for any function $f$ which is locally conformal throughout the Riemann sphere, except for an at most countable set of points in which $f'(z) = 0$, or which are multiple poles of $f$, or which are isolated essential singularities.

In the case of functions having essential singularities, the Big Picard Theorem has been instrumental in finding fundamental domains [2]. Once this achieved, by using the technique of simultaneous continuations [7], we obtained an essentially enriched type of Picard Theorem, saying that *every neighborhood V of an essential singularity of an analytic function f contains infinitely many fundamental domains*, i.e., after Ahlfors definition, *domains which are conformally mapped by f onto the whole complex plane with a slit.*

Obviously, *f* takes in *V* every value, except possibly those at the ends of the slit (lacunary values) infinitely many times. Moreover, we have infinitely many disjoint domains in *V* which are mapped conformally by *f* onto a neighborhood of every such value.

The Riemann Zeta function has a unique essential singularity at the point ∞ of the Riemann sphere. Thus, its fundamental domains should accumulate to infinity and only there, in the sense that every compact set from the complex plane intersects only a finite number of fundamental domains of this function, while in the exterior of a compact set there are infinitely many such domains.

## 2. The Pre-Image by $\zeta$ of the Real Axis

In order to make the paper self contained, let us repeat some of the results obtained in [3]. The pre-image by $\zeta$ of the real axis can be viewed as simultaneous continuation over the real axis from a real value starting from all the points in which that value is assumed.

By the Big Picard Theorem, every value $z_0$ from the *z*-plane ($z = \zeta(s)$), if it is not a lacunary value, is taken by the function $\zeta$ in infinitely many points $s_n$ accumulating to ∞ and only there. This is true, in particular, for $z_0 = 0$. By [9], the Zeta function has only simple zeros, hence a small interval *I* of the real axis containing 0 will have as pre-image by $\zeta$ the union of infinitely many Jordan arcs $\gamma_n$ passing each one through a zero $s_n$ of $\zeta$, and vice-versa, every zero $s_n$ belongs to some arcs $\gamma_n$. Since $\zeta(\sigma) \in \mathbb{R}$, for $\sigma \in \mathbb{R}$, and by the formula (5), the *trivial zeros* of $\zeta$ are simple zeros and the arcs $\gamma_n$ corresponding to these zeros are intervals of the real axis, if *I* is small enough.

Due to the fact that $\zeta$ is analytic (except at $s = 1$), between two consecutive trivial zeros of $\zeta$ there is at least one zero of the derivative $\zeta'$, i.e. at least one branch point of $\zeta$. Since we have also $\zeta'(\sigma) \in \mathbb{R}$ for $\sigma \in \mathbb{R}$, if we perform simultaneous continuations over the real axis of the components included in $\mathbb{R}$ of the pre-image by $\zeta$ of *I*, we will encounter at some moments those branch points and the continuations follow on unbounded curves crossing the real axis at the respective points. As shown in [3], the respective curves are unbounded and they do not intersect each other.

Only the continuation of the interval containing the zero $s = -2$ stops at the unique pole $s = 1$, since $\lim_{\sigma \nearrow 1} \zeta(\sigma) = \infty$. Similarly, if instead of $z_0 = 0$ we take another real $z_0$ greater than 1 and perform the same operations, since $\lim_{\sigma \searrow 1} \zeta(\sigma) = \infty$, the continuation over the interval $(1, \infty)$ stops again at $s = 1$. In particular, the pre-image by $\zeta$ of this interval can contain no zero of $\zeta$.

Thus, if we color for example red, the pre-image by $\zeta$ of the negative real half axis and let black the pre-image of the positive real half axis, then all the components of the pre-image of the interval $(1, +\infty)$ will be black, while those of the interval $(-\infty, 1)$ will have a part red and another part black, the junction of the two colors corresponding to a zero of $\zeta$ (trivial or not).

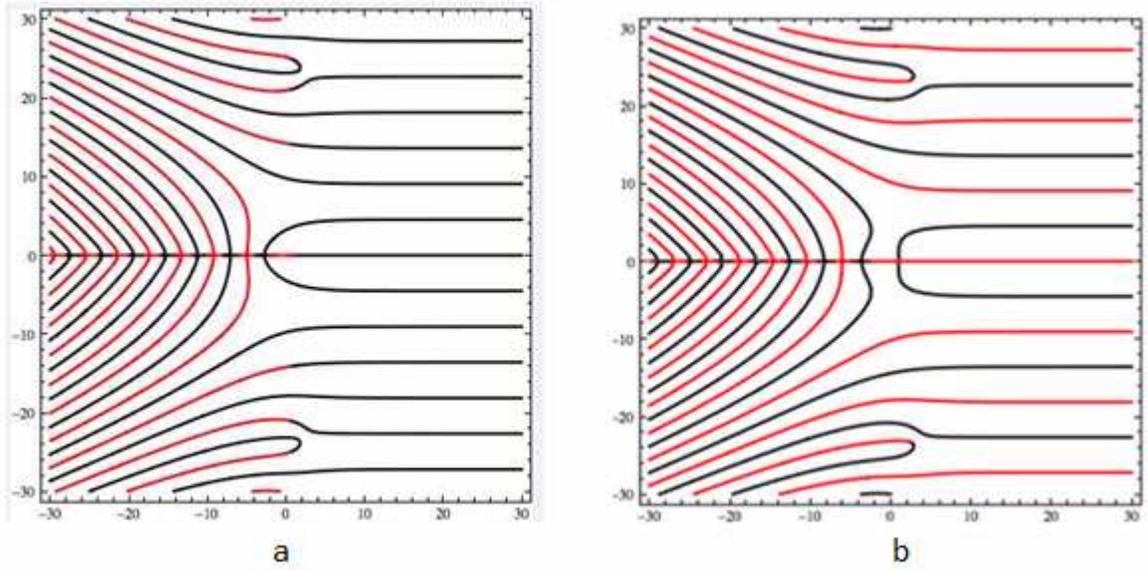

Fig. 1

Fig.1 represents the pre-images by $\zeta$ (Fig. 1a) and by $\zeta'$ (Fig. 1b) of the real axis in boxes of $[-30,30] \times [-30,30]$.

We notice the existence of branch points on the real axis and their color alternation, as well as the trivial zeros between them. If we superpose the two pictures, the branch points of the first, will coincide with the zeros of the second.

The components passing through non trivial zeros of the pre-image of the real axis form a more complex configuration, which has a lot to do with the special status of the value $z = 1$: on one hand this value is taken in infinitely many points of the argument plane, and on the other hand it behaves like a lacunary value, since it is obtained as a limit as $s$ tends to infinity on some unbounded curves. We called it *quasi lacunary value*.

For the function $\zeta'$ the value $z = 0$ is quasi lacunary.

Due to the symmetry with respect to the real axis ($\zeta(\bar{s}) = \overline{\zeta(s)}$), it is sufficient to deal only with the upper half plane. Let $x_0 \in (1, +\infty)$ and let $s_k \in \zeta^{-1}(\{x_0\}) \setminus \mathbb{R}$. Continuation over $(1, +\infty)$ from $s_k$ is either an unbounded curve $\Gamma'_k$ such that $\lim_{\sigma \to +\infty} \zeta(\sigma + it) = 1$, by (1), and $\lim_{\sigma \to -\infty} \zeta(\sigma + it) = +\infty$, where $s = \sigma + it \in \Gamma'_k$, or there are points $u$ such that $\zeta(u) = 1$, thus the continuation can take place over the whole real axis.

**Theorem 1**. *Consecutive curves $\Gamma'_k$ and $\Gamma'_{k+1}$ form strips $S_k$ which are infinite in both directions. The function $\zeta$ maps these strips (not necessarily bijectively) onto the complex plane with a slit alongside the interval $[1, +\infty)$ of the real axis.*

*Proof:* Here $k \in \mathbb{N}$ is chosen in such a way that if $it_k \in \Gamma'_k$, then the sequence $(t_k)$ is increasing. If two consecutive curves $\Gamma'_k$ and $\Gamma'_{k+1}$ met at a point $s$, one of the domains bounded by them would be mapped by $\zeta$ onto the complex plane with a slit alongside the real axis from 1 to $\zeta(s)$. Such a domain must contain a pole of $\zeta$, but this cannot happen, since the only pole of $\zeta$ is $s = 1$. When a point $s$ travels on $\Gamma'_k$ and then on $\Gamma'_{k+1}$ leaving the strip at left, $\zeta(s)$ moves on the real axis from 1 to $\infty$ and back. $\square$

When the continuation can take place over the whole real axis, we obtain unbounded curves each containing a non trivial zero of $\zeta$ and a point $u$ with $\zeta(u) = 1$. Such a point $u$ is necessarily interior to a strip $S_k$ since the border of every $S_k$, which is included in $\zeta^{-1}((1,+\infty))$, and the set $\zeta^{-1}(\{1\})$ are disjoint.

Let us denote by $u_{k,j}$ the points of $S_k$ for which $\zeta(u_{k,j}) = 1$, by $\Gamma_{k,j}, j \neq 0$ the components of $\zeta^{-1}(\mathbb{R})$ containing $u_{k,j}$ and by $s_{k,j}$ the non trivial zero of $\zeta$ situated on $\Gamma_{k,j}$. As shown in [3] Theorem 7, every $S_k$ contains a unique component $\Gamma_{k,0}$ with the property that $\lim_{\sigma \to +\infty} \zeta(\sigma + it) = 1$ and $\lim_{\sigma \to -\infty} \zeta(\sigma + it) = -\infty$. where $\sigma + it \in \Gamma_{k,0}$ (i.e. which is projected bijectively by $\zeta$ onto the interval $(-\infty, 1)$) and a finite number $j_k - 1$ of components such that each one is projected bijectively by $\zeta$ onto the whole real axis. Therefore $S_k$ contains $j_k$ curves $\Gamma_{k,j}$, $j_k$ non trivial zeros and $j_k - 1$ points $u_{k,j}$ with $\zeta(u_{k,j}) = 1$. We call $S_k$ a $j_k$-strip. Here $j \in \mathbb{Z}$ is chosen in such a way that $\Gamma_{k,j}$ and $\Gamma_{k,j+1}$ are consecutive in the same sense as $\Gamma'_k$.

**Theorem 2**: *When the continuation takes place over the whole real axis, the components $\Gamma_{k,j}$ are such that the branches corresponding to both the positive and the negative real half axis contain only points $\sigma + it$ with $\sigma < 0$ for $|\sigma|$ big enough.*

*Proof:* A point traveling in the same direction on a circle $\gamma$ centered at the origin of the $z$-plane meets consecutively the positive and the negative real half axis. Thus the pre-image of $\gamma$ should meet consecutively the branches corresponding to the pre-image of the positive and the negative real half axis, which is possible only if the condition of the theorem is fulfilled.$\square$

In what follows, this *color alternation condition* will appear repeatedly. One of its immediate applications is that the real zeros of $\zeta'$ are simple and they alternate with the trivial zeros of $\zeta$. Indeed, since the color change can happen only at a zero of $\zeta$ and at $s = 1$ and by the formula (5) the trivial zeros of $\zeta$ are simple, between two of them there must be one and only one zero of $\zeta'$, otherwise the color alternation condition would be violated.

Another application is that those components of pre-images of circles centered at the origin of radius $\rho > 1$ which cross a $\Gamma'_k$, will continue to cross alternatively red and black components of the pre-image of the real axis indefinitely, i.e. they are unbounded components.

**Theorem 3**. *Every strip $S_k$ contains a unique unbounded component of the pre-image of the unit disc.*

*Proof:* To see this, it is enough to take the pre-image of a ray making an angle $\alpha$ with the real half axis and let $\alpha \to 0$. A point $s \in \zeta^{-1}(\{z\})$ with $|z| = 1$, $\arg z = \alpha$, tends to $\infty$ as $\alpha \to 0$ if and only if the corresponding component of the pre-image of that ray tends to $\Gamma_{k,0}$ as $\alpha \to 0$, which happens if and only if the component of the pre-image of the closed unit disc containing the point $s$ is unbounded. The uniqueness of $\Gamma_{k,0}$ implies the uniqueness of such a component.$\square$

An even better way to visualize the geometry of conformal mappings realized by $\zeta$ is to take the pre-image of an orthogonal net formed with rays passing through the origin and circles centered at the origin. The annuli so formed are colored as in Fig. 2, where the points in the annuli and their pre-images by $\zeta$ have the same color, saturation and brightness. The color saturation in the annuli is decreasing clockwise. We show in the Fig. 2 (a,b,c) annuli of radii up to $10^4$, although we investigated strips around the critical line with values of $t$ up to $10^9$ (see pictures at the end ). We notice that to the quadrilaterals of the net from the $z$-plane correspond quadrilaterals in the argument plane of $\zeta$ which not only have the same colors but also the same conformal module. By using a convenient scale on the rays, we can arrange that all these quadrilaterals have the same module. The components of the pre-images of the four quadrilaterals having a vertex in $z = 1$ are unbounded if they don't have an $u_{k,j}$ as a vertex. It makes sense to say that the module of such an unbounded quadrilateral is that of its (bounded) image by $\zeta$.

The theorems above guarantee that the landscape shown in Fig. 2(d,e,f) repeats itself indefinitely with variations regarding only the numbers of zeros in a strip $S_k$.

On the other hand, this landscape gives one an idea how difficult can it be to tackle the Riemann Hypothesis by number theoretical methods and suggests to rather use conformal mappings tools.

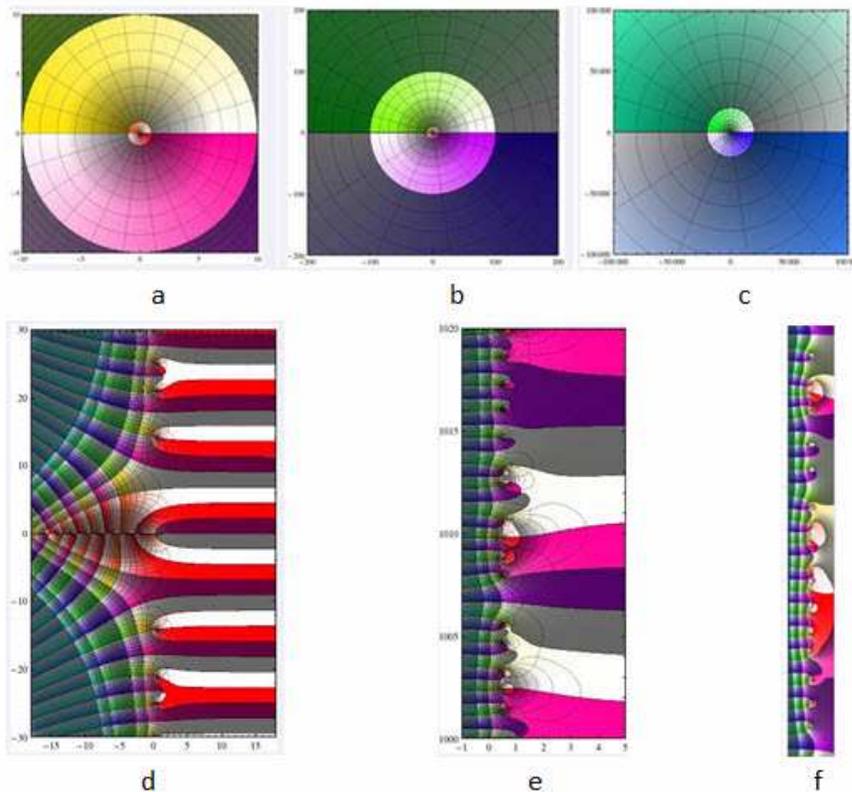

Fig. 2

3. **Fundamental Domains of the Riemann Zeta Function**

If $S_k$ is a $j_k$-strip, then the pre-image of a small circle $\gamma$ of radius $\rho$ centered at the origin will have $j_k$ components situated in $S_k$ which are closed Jordan curves each containing a unique zero of $\zeta$. If $\rho$ is small enough, then these curves are disjoint. As $\rho$ increases, the curves expand. When two of them touch each other at a points $v_{k,j}$, they fuse into a unique closed curve containing the two respective zeros. Continuing to increase $\rho$, some other points $v_{k,j}$ are reached. These are branch points of $\zeta$. It has been shown in [3] and [9] that there are exactly $j_k - 1$ points $v_{k,j}$ in every $j_k$-strip. The intersection of $S_k$ with the pre-image of the segment between 1 and $\zeta(v_{k,j})$ is a set of unbounded curves starting at $u_{k,j}$ or arcs connecting two such points. These curves and arcs together with the pre-image of the interval $[1,+\infty)$ of the real axis bound sub-strips included in $S_k$ which are fundamental domains $\Omega_{k,j}$ of $\zeta$. Every fundamental domain $\Omega_{k,j}$ contains a unique simple zero $s_{k,j}$ of $\zeta$.

The simplicity of the zeros of $\zeta$ has been proved in [9] by using the so-called *intertwined curves*. Namely, it has been shown that the components $\Upsilon'_k$ of the pre-image by $\zeta'$ of the interval $(-\infty, 0)$ of the real axis such that $\lim_{\sigma \to -\infty} \zeta'(\sigma + it) = -\infty$ and $\lim_{\sigma \to +\infty} \zeta'(\sigma + it) = 0$, $\sigma + it \in \Upsilon'_k$ form infinite strips $\Sigma_k$ containing $j_k - 1$ simple zeros of $\zeta'$. These zeros belong to the components $\Upsilon_{k,j}$, $j \neq 0$ of the pre-image by $\zeta'$ of the real axis. There is also in $\Sigma_k$ a unique component $\Upsilon_{k,0}$ of the pre-image of the interval $(0,+\infty)$. These components follow also the color alternation rule. It has been proved in [9] that there is a one to one correspondence between those $\Gamma'_k$ and $\Upsilon'_k$, respectively $\Gamma_{k,j}$ and $\Upsilon_{k,j}$ which intersect each other (intertwined curves). Also, by using different colors for the components of the pre-image by $\zeta$ and by $\zeta'$ of the positive and negative real half axes, a *color matching rule* has been proved in [9]. The simplicity of the zeros of Zeta function, as well as of those of any derivative of $\zeta$ is a consequence of the two rules: the color alternation rule for $\zeta$ and $\zeta'$ and the color matching rule.

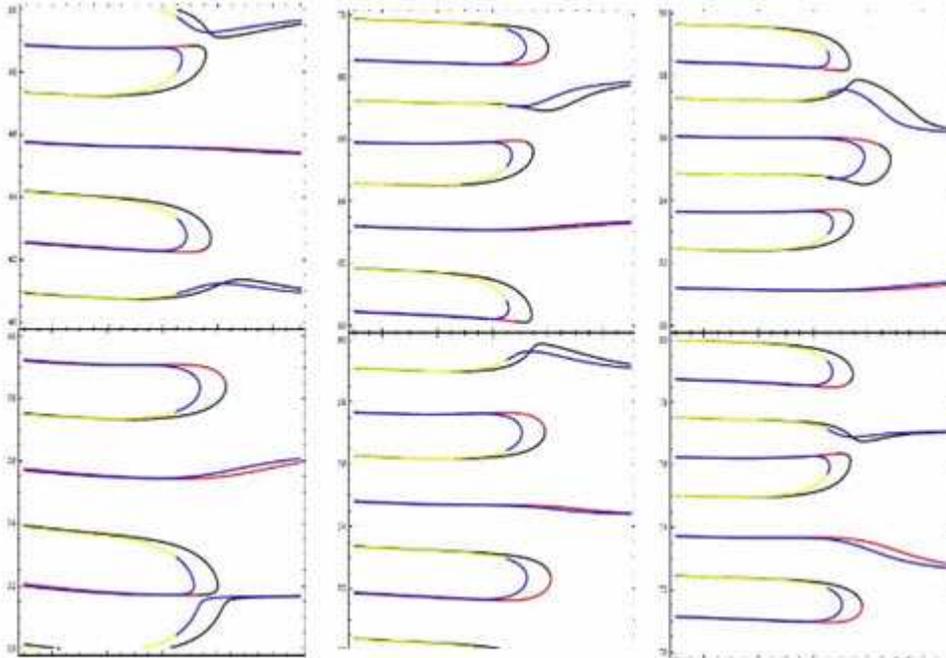

Fig. 3

Fig. 3 above illustrates the concept of intertwined curves, parts of which are shown in a box [-10,10]x[30,90].

The color matching rule means that *blue matches red and yellow matches black* for all the components except for $\Gamma_{k,0}$ and $\Upsilon_{k,0}$.

For the seek of space economy, the height of the box has been divided in three.

An inspection of the curves $\Gamma_{k,j}$ from the Fig. 3 shows that they have all the turning point

(the point in which the tangent to $\Gamma_{k,j}$ is vertical) on the blue part, i.e. on the part corresponding to the positive real half axis. We will prove next that this happens for any curve $\Gamma_{k,j}$.

**Theorem 4**. *Let $x \to s_{k,j}(x) = \sigma_{k,j}(x) + it_{k,j}(x)$ be the parametric equation of $\Gamma_{k,j}$ such that*

$$\zeta(s_{k,j}(x)) = x, \ x \in \mathbb{R}, \ s_{k,j} = \sigma_{k,j} + it_{k,j} = s_{k,j}(0),$$

*If $\sigma'_{k,j}(x_0) = 0$, then $x_0 > 0$*

*Proof:* In terms of the zeros $s_{k,j}$, the *Hadamard Product Formula* (see Wikipedia / Riemann Zeta Function) can be written as:

$$\zeta(s) = [1 / 2(s-1)\Gamma(1+s/2)]e^{As} \prod_{k=1}^{\infty} \prod_{j \in J_k} (1 - s/\bar{s}_{k,j})(1 - s/s_{k,j}) \exp\{s(1/\bar{s}_{k,j} + 1/s_{k,j})\},$$

where $A = 0.549269234\ldots$
Let

$$\varphi_{k,j}(x) = [1 - s_{k,j}(x)/\bar{s}_{k,j}][1 - s_{k,j}(x)/s_{k,j}]$$

Since all the zeros of $\zeta$ are simple, in the neighborhood of $s_{k,j}$, the dominating factor of $\zeta(s_{k,j}(x))$ is $\varphi_{k,j}(x)$. We have

$$\varphi_{k,j}(x) = 1 - 2s_{k,j}(x)\sigma_{k,j}/|s_{k,j}|^2 + s_{k,j}^2(x)/|s_{k,j}|^2$$
$$= [1/|s_{k,j}|^2][(\sigma_{k,j}(x) - \sigma_{k,j})^2 + t_{k,j}^2 - t_{k,j}^2(x) + 2t_{k,j}(x)(\sigma_{k,j}(x) - \sigma_{k,j})i].$$

Thus, if $t_{k,j}(x) < t_{k,j}$, then $\operatorname{Re} \varphi_{k,j}(x) > 0$ and for $|\sigma_{k,j}(x) - \sigma_{k,j}|$ small enough, if $t_{k,j}(x) > t_{k,j}$ then $\operatorname{Re} \varphi_{k,j}(x) < 0$. Also $\operatorname{Im} \varphi_{k,j}(x) > 0$ if and only if $\sigma_{k,j}(x) > \sigma_{k,j}$.

Suppose that for $x < x_0 < 0$, we have $\sigma_{k,j}(x) > \sigma_{k,j}$, $t_{k,j}(x) < t_{k,j}$, i.e. the curve $\Gamma_{k,j}$ turns back at $s(x_0)$, where $x_0 < 0$, i.e. on the branch corresponding to the negative half axis, which is situated below the other branch. Then

$$\pi/2 < \arg s'_{k,j}(0) < \pi, \quad 0 < \arg \varphi_{k,j}(x) < \pi/2 \qquad (*)$$

Since

$$\varphi'_{k,j}(x) = \frac{2}{|s_{k,j}|^2}[s_{k,j}(x) - \sigma_{k,j}]s'_{k,j}(x) \text{ and } s_{k,j}(0) = s_{k,j}, \text{ we have}$$

$$\varphi'_{k,j}(0) = \frac{2}{|s_{k,j}|^2}[s_{k,j} - \sigma_{k,j}]s'_{k,j}(0) = \frac{2it_{k,j}}{|s_{k,j}|^2}s'_{k,j}(0),$$

hence

$$\arg \varphi'_{k,j}(0) = \frac{\pi}{2} + \arg s'_{k,j}(0) \qquad (**)$$

Since $\varphi'_{k,j}(0) = \lim_{x \to 0}[\varphi_{k,j}(0) - \varphi_{k,j}(x)]/x = \lim_{x \to 0}[-\varphi_{k,j}(x)]/x$ and $x < 0$, we have that

$$0 \leq \arg \varphi'_{k,j}(0) \leq \pi/2$$

which does not agree with $(*)$ and $(**)$.

If the curve turns back on the part corresponding to the positive half axis, and the relative position of the two branches is the same, then for $x > 0$ small enough, we have $\sigma(x) > \sigma_{k,j}$, $t(x) > t_{k,j}$, hence $\pi/2 < \arg \varphi_{k,j}(x) < \pi$ and $0 < \arg s'(0) < \pi/2$.

Then $\pi/2 < \arg \varphi'_{k,j}(0) < \pi$, which agrees with the formula $(**)$. Similar arguments are valid when the relative position of two branches is reversed. The conclusion is that all the curves $\Gamma_{k,j}$ turn back at points corresponding to $x > 0$.

### 4. The Proof of Riemann Hypothesis

From the relation (5) and the fact that $\zeta(\bar{s}) = \overline{\zeta(s)}$, it can be easily inferred that the non trivial zeros of $\zeta$ appear in quadruplets: $s, \bar{s}, 1 - s$ and $1 - \bar{s}$. Proving the RH is equivalent to showing that $\zeta$ cannot have two zeros of the form $s_1 = \sigma_1 + it$, $s_2 = \sigma_2 + it$ in the critical strip. In particular, if we suppose that $s_1 = \sigma + it$ and $s_2 = 1 - \bar{s}_1 = 1 - \sigma + it$, $0 < \sigma < 1/2$ are zeros of $\zeta$, then $s_1 = s_2$, i.e. $\sigma = 1 - \sigma$, hence $\sigma = 1/2$. In other words, all the non trivial zeros of $\zeta$ are of the form $s = 1/2 + it$.

There are two situations to be examined, namely a), when $s_1$ and $s_2$ are consecutive zeros in the same strip $S_k$ and b), when $s_1$ is the last zero of $S_k$ and $s_2$ is the first zero of $S_{k+1}$.

Once assured that $\sigma = 1/2$ in the two cases, a simple induction argument guarantees that this happens for all non trivial zeros in the upper half plane, and then, due to the symmetry $\zeta(\bar{s}) = \overline{\zeta(s)}$, the truth of the Riemann Hypothesis follows.

For the seek of completeness, let us present that induction argument.

Since $S_1$ contains only one zero $s_{1,0}$, the case a) is void for it.

Then, applying the case b) to $S_1$ and $S_2$ one concludes that $s_{1,0} = 1/2 + it_{1,0}$, for some positive number $t_{1,0}$.

Suppose that for an arbitrary $k$ all the zeros in $S_k$ are of the form $s_{k,j} = 1/2 + it_{k,j}$ and apply the case b) to the last zero $s_{k,j'}$ from $S_k$ and the first zero $s_{k+1,j''}$ from $S_{k+1}$.

The conclusion is that $s_{k+1,j''} = 1/2 + it_{k+1,j''}$ for some positive number $t_{k+1,j''}$.

Then, applying the case a) to every couple of consecutive zeros from $S_{k+1}$ one can draw the conclusion that all the zeros in $S_{k+1}$ are of the form $1/2 + it_{k+1,j}$ for some positive numbers $t_{k+1,j}$, which completes the induction.

*Case a).*

Suppose that the respective zeros are $s_{k,j-1} = \sigma_{k,j-1} + it$ and $s_{k,j} = \sigma_{k,j} + it$. Let $\Gamma_{k,j}$ and $\Upsilon_{k,j}$, respectively $\Gamma_{k,j-1}$ and $\Upsilon_{k,j-1}$ be the corresponding intertwined curves with $s_{k,j} \in \Gamma_{k,j}$ and $v_{k,j} \in \Upsilon_{k,j}$, respectively $s_{k,j-1} \in \Gamma_{k,j-1}$ and $v_{k,j-1} \in \Upsilon_{k,j-1}$.

Fig. 4 below represents such a hypothetical situation.

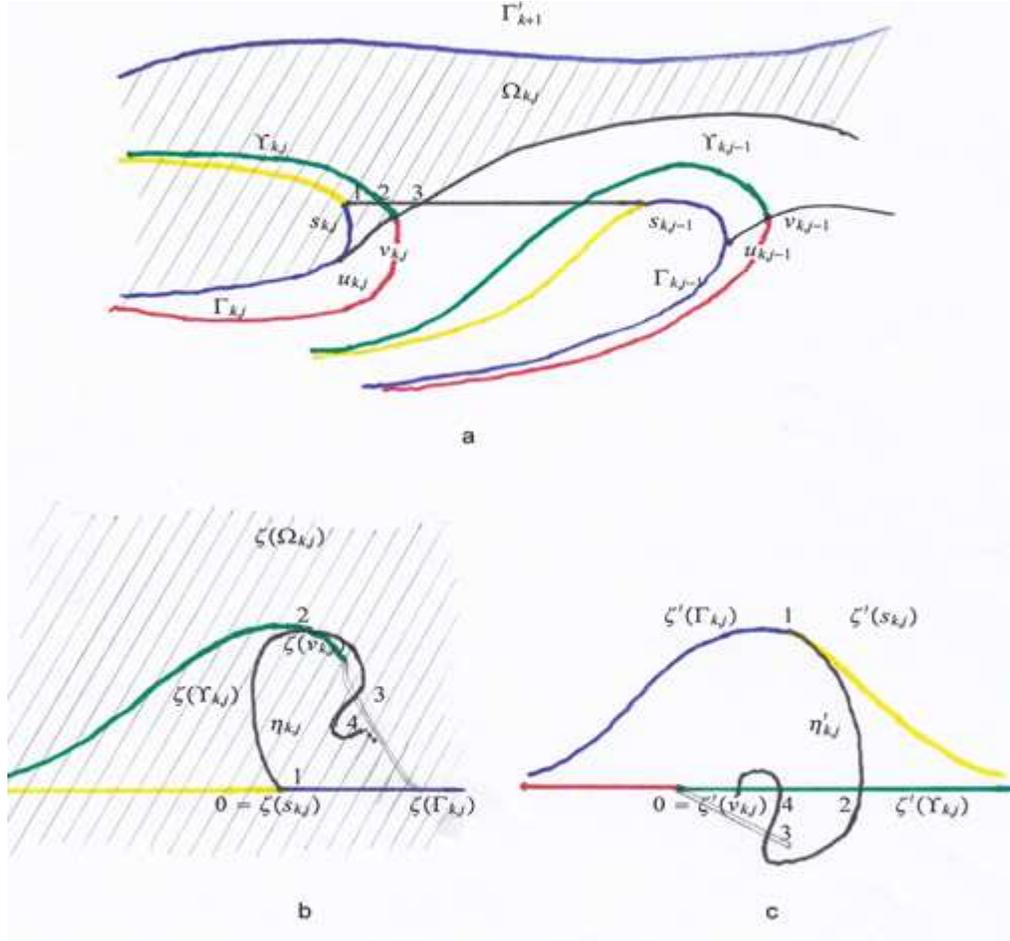

Fig. 4

If $\lambda \rightarrow (1 - \lambda)s_{k,j-1} + \lambda s_{k,j}, 0 \leq \lambda \leq 1$ is the parametric equation of the segment $I$ between $s_{k,j}$ and $s_{k,j-1}$, then

$$z(\lambda) = \zeta((1 - \lambda)s_{k,j} + \lambda s_{k,j-1}), 0 \leq \lambda \leq 1 \quad (6)$$

is a parametric equation of a closed curve $\eta_{k,j}$ passing through the origin. The curve is obviously smooth except possibly at the origin and at $v_{k,j}$, if $v_{k,j} \in I$. Since $v_{k,j}$ is a simple zero of $\zeta'$ we have in $v_{k,j}$ a *star configuration,* (as in [1], page 133) of two orthogonal curves, the image by $\zeta$ of one of which passes through the origin. The tangent to this curve $\eta_{k,j}$ at $v_{k,j}$ still exists. In any case we can differentiate in (6) with respect to $\lambda$, $0 < \lambda < 1$ and we get:

$$z'(\lambda) = \zeta'((1-\lambda)s_{k,j-1} + \lambda s_{k,j})(s_{k,j-1} - s_{k,j}). \qquad (7)$$

Since $s_{k,j-1} - s_{k,j} = \sigma_{k,j-1} - \sigma_{k,j} > 0$, we have

$$\arg z'(\lambda) = \arg \zeta'((1-\lambda)s_{k,j} + \lambda s_{k,j-1}) \qquad (8)$$

which shows that the tangent to $\eta_{k,j}$ at $z(\lambda)$ makes an angle with the positive real half axis which is equal to the argument of $\zeta'((1-\lambda)s_{k,j} + \lambda s_{k,j-1})$. This last point describes an arc $\eta'_{k,j}$ starting at $\zeta'(s_{k,j})$ and ending at $\zeta'(s_{k,j-1})$, when $\lambda$ varies from 0 to 1. Due to the color alternation rule and the color matching rule, both ends of $\eta'_{k,j}$ must be situated in the upper half plane. The right hand side in (8) exists also at the ends of this arc, which means that the limits as $\lambda \searrow 0$ and as $\lambda \nearrow 1$ of the left hand side also exist. These limits are two vectors starting at the origin with the first pointing to the upper half plane and the second pointing to the lower half plane. Indeed, the oriented interval $I$ exits from the *parabola-like* curve $\Gamma_{k,j}$ and enters the *parabola-like* curve $\Gamma_{k,j-1}$ whose interiors are conformally mapped by $\zeta$ onto the lower half plane. This contradicts the equality (8) and the fact that both ends of $\eta'_{k,j}$ are situated in the upper half plane.

The conclusion is that a configuration as postulated in Fig. 4a is impossible.

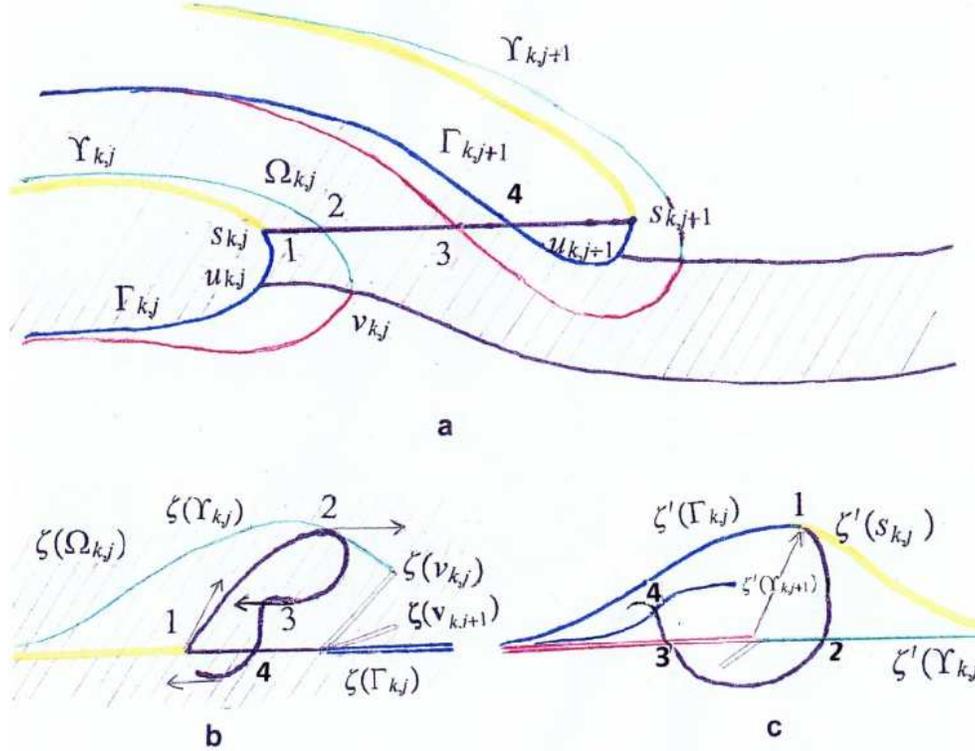

Fig 4'

Let us deal with a hypothetical situation as represented in Fig. $4'$a. The points 1,2,3,... are mapped in tandem by $\zeta$ and $\zeta'$ into the points represented by the same numbers in Fig. $4'$b, respectively $4'$c. The number 1 is the origin in Fig. $4'$b and it is $\zeta'(s_{k,j})$ in Fig $4'$c. Point 2 is the intersection with the positive half axis of $\eta'_{k,j}$, hence $\arg \zeta'((1+\lambda)s_{k,j} + \lambda s_{k,j-1}) = 0$ at the point 2, thus $\arg z'(\lambda) = 0$, which means that the tangent to $\eta_{k,j}$ at that point is horizontal and positively oriented with respect to the real axis. Point 3 corresponds in Fig. $4'$a to the intersection of $I$ with $\Upsilon_{k,j+1}$ and is situated in Fig. $4'$c on the negative real half axis. Then $\eta_{k,j}$ has a horizontal tangent at the corresponding point 3 oriented in the negative direction of the real axis. At point 4 the interval $I$ crosses $\Gamma_{k,j+1}$, hence $\eta'_{k,j}$ will cross the curve $\zeta'(\Gamma_{k,j+1})$, while $\eta_{k,j}$ will cross the real axis. In order for $\eta_{k,j}$ to reach the origin, it must have another point with horizontal tangent. This forces $\eta'_{k,j}$ to turn back to the real axis, crossing in its way again $\zeta'(\Gamma_{k,j+1})$ and so on. This means that $I$ intersects again $\Gamma_{k,j+1}$ and $\Upsilon_{k,j+1}$, before reaching $s_{k,j+1}$, which is absurd. Thus a configuration like that postulated in Fig. $4'$a is impossible.

The final conclusion is that $\zeta$ cannot have two zeros of the form $s_{k,j} = \sigma + it_{k,j}$ and $s_{k,j+1} = 1 - \sigma + it_{k,j+1}$.

*Case b).*

Suppose now that we deal with the zeros $s_{k,j'}$ and $s_{k+1,j}$ belonging to the adjacent strips $S_k$, respectively $S_{k+1}$. We make similar notations as in the case a) and we seek to find a contradiction by inspecting the images by $\zeta$ and by $\zeta'$ of the interval $I$ between the two zeros whose parametric equation is $\lambda \to (1-\lambda)s_{k+1,j} + \lambda s_{k,j'}$.

Fig. 5 below illustrates the new hypothetical situation.

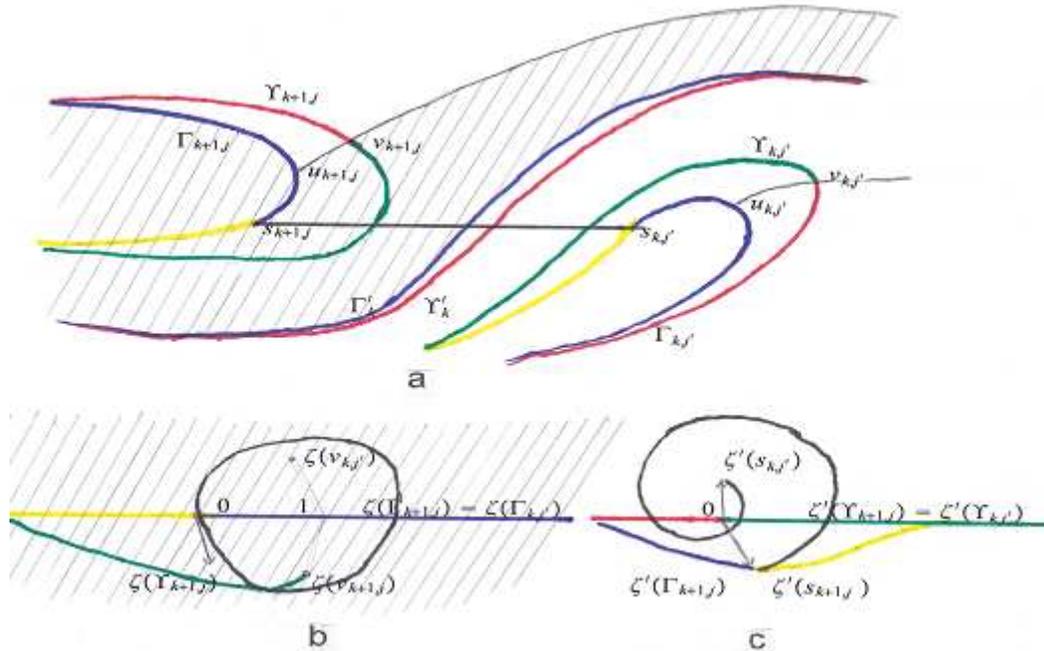

Fig. 5

This time $\zeta'(s_{k,j'})$ belongs to the upper half plane, while $\zeta'(s_{k+1,j})$ belongs to the lower half plane. We notice again that $\lim_{\lambda \searrow 0} z'(\lambda)$ and $\lim_{\lambda \nearrow 1} z'(\lambda)$ exist and represent the half-tangents to $\eta_{k,j}$ at the origin. Both of these vectors point to the lower half plane, which contradicts the fact that $\zeta'(s_{k,j'})$ belongs to the upper half plane.

Thus, a configuration like that in Fig. 5a is impossible.

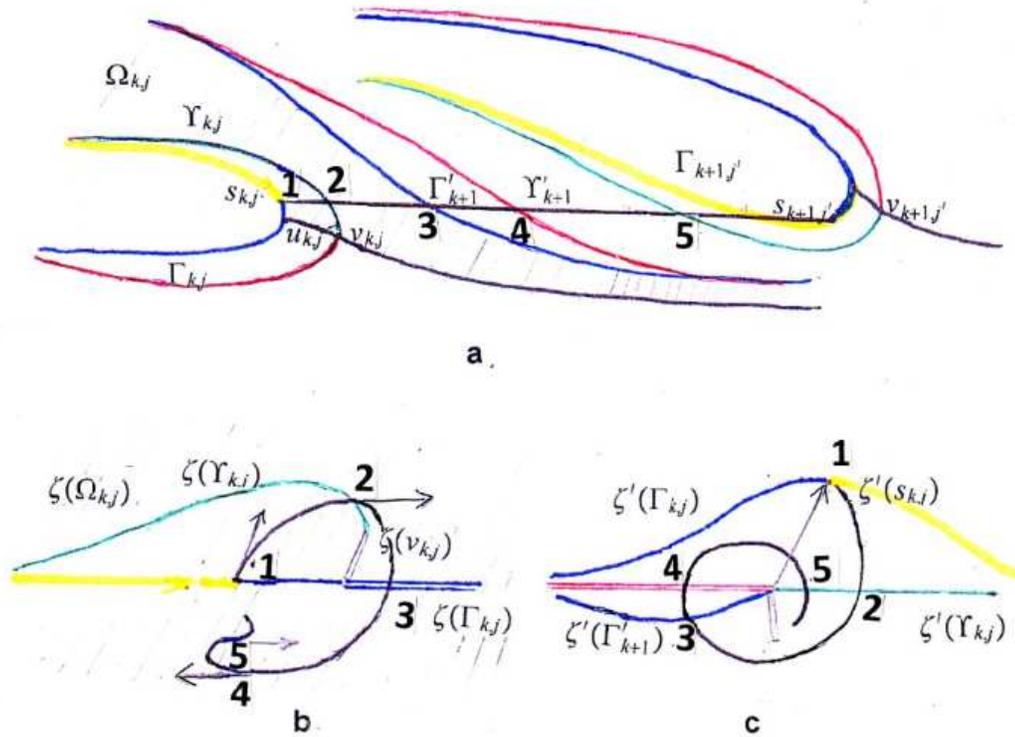

Fig. 5'

The argument corresponding to the configuration shown in Fig. 5' goes as follows. The position vectors of the points on $\eta'_{k,j}$ point to the lower half plane after point 5. However, $\eta_{k,j}$ must turn to the origin after the corresponding position 5, which forces the tangent to it to point towards the upper half plane, contradicting the equality (8).

Thus, a configuration like that shown in Fig. 5'a is not possible.

Consequently, all the imaginable positions of two consecutive non trivial zeros having the same imaginary part bring us to contradictions.

This completes the proof of Riemann Hypothesis □.

**Acknowledgments**: The author is grateful to Cristina Ballantine and Florin Alan Muscutar for providing computer generated graphics.

**References**

[1] Ahlfors, L.V., *Complex Analysis*, International Series in Pure and Applied Mathematics (1979)

[2] Andreian Cazacu, C. and Ghisa, D., *Global Mapping Properties of Analytic Functions*, Proceedings of the 7-th ISAAC Congress, London, U.K., 3-12, (2009)

[3] Andreian Cazacu, C. and Ghisa, D., *Fundamental Domains of Gamma and Zeta Functions*, International Journal of Mathematics and Mathematical Sciences, 985323, 21 pages, Volume 2011

[4] Ballantine, C. and Ghisa, D., *Color Visualization of Blaschke Product Mappings*, Complex Variables and Elliptic Equations, Vol. 55, Nos. 1-3, 201-217, (2010)

[5] Ballantine, C. and Ghisa, D., *Color Visualization of Blaschke Self-Mappings of the Real Projective Plan*, Proceedings of the International Conference on Complex Analysis and Related Topics, Alba Iulia, Romania (2008), Rev. Roum. Math. Pure Appl., Vol. 54, Nos. 5-6, 375-394, (2009)

[6] Ballantine, C. and Ghisa, D., *Global Mapping Properties of Rational Functions*, Proceedings of the 7-th ISAAC Congress, London, U.K., 13-22, (2009)

[7] Barza, I. and Ghisa, D., *The Geometry of Blaschke Product Mappings*, Proceedings of the 6-th International ISAAC Congress, Ankara (2007) in Further Progress in Analysis, Editors H.G.W.Begehr, A.O. Celebi & R.P.Gilbert, World Scientific, 197–207 (2008)

[8] Derbyshire, J., *Prime Obsession, Bernhard Riemann and the Greatest Unsolved Problem in Mathematics*, Plume, 2003

[9] Ghisa, D. and F.A.Muscutar, *The Riemann Zeta Function Has Only Simple Zeros,* in Proceedings of the XVIth Conference on Analytic Functions and Related Topics, Chelm, Poland, 2011

**Addendum**

We add some more pictures, in order to illustrate some facts and to show some surprising phenomena.

From these pictures it appears that the number $j_k$ of fundamental domains included in $S_k$ increases (on the average) with $k$. These numbers are, for $k$ between 1 and 11 respectively: 1, 2, 2, 3, 3, 3, 4, , 3, 4, 4, 4 (see Fig. 6). The strip $S_k$ corresponding to $t$ in the range of 1000 has seven fundamental domains (Fig.2 e), that corresponding to $t$ in the range of $10^6$ has 17 fundamenral domains (fig.2 f) and that corresponding to $t$ in the range of $10^9$ has 29 fundamental domains (Fig. 7).

In the strips in this range two unexpected phenomena can be seen, namely a bounded component of the pre-image of the unit disc situated in front of the unbounded component in the interval for $t$ between $10^9 + 3$ and $10^9 + 4$ as well as curves $\Gamma_{k,j}$ embracing some other such curves, as in the intervals for $t$ between $10^9 + 11$ and $10^9 + 13$, respectively $10^9 + 18$ and $10^9 + 19$. These phenomena do not contradict the theorems proven in [3] and [9].

We notice also that the width of $S_k$ in the range of the critical strip is approximately the same for every $k$, namely in the range of 10, which means that the distance between two consecutive non trivial zeros tends (in average) to zero as $t \to \infty$, yet it can never be zero.

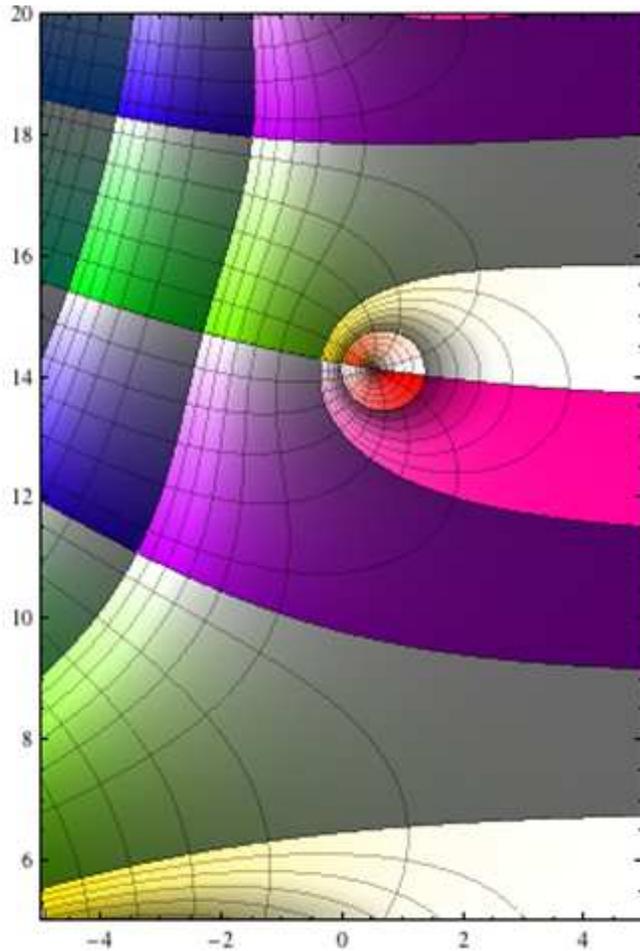

This is the only 1-strip in the upper half plane.

It is a *fundamental domain*, mapped by Zeta function onto the z-plane with a slit alongside the real axis from 1 to infinity.

An arbitrary *m-strip* contains *m fundamental domains* and for every one a *supplementary slit* in the unit disc is added.

Fig.6 a

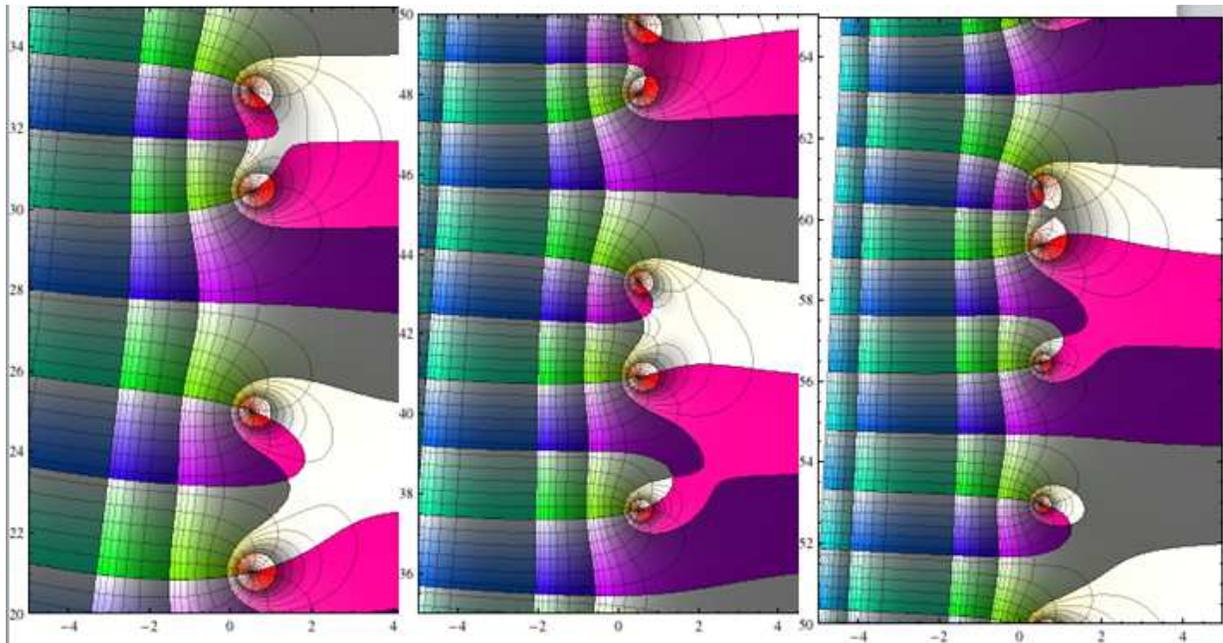

Fig. 6 b

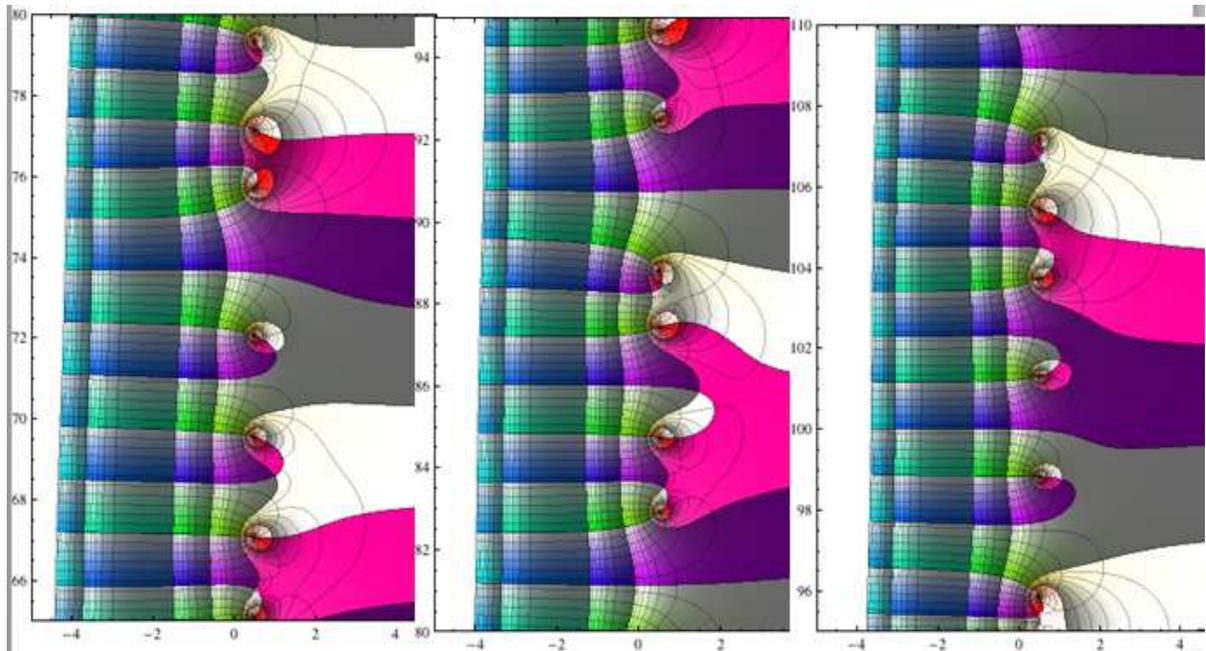

Fig.. 6 c

Fig. 7

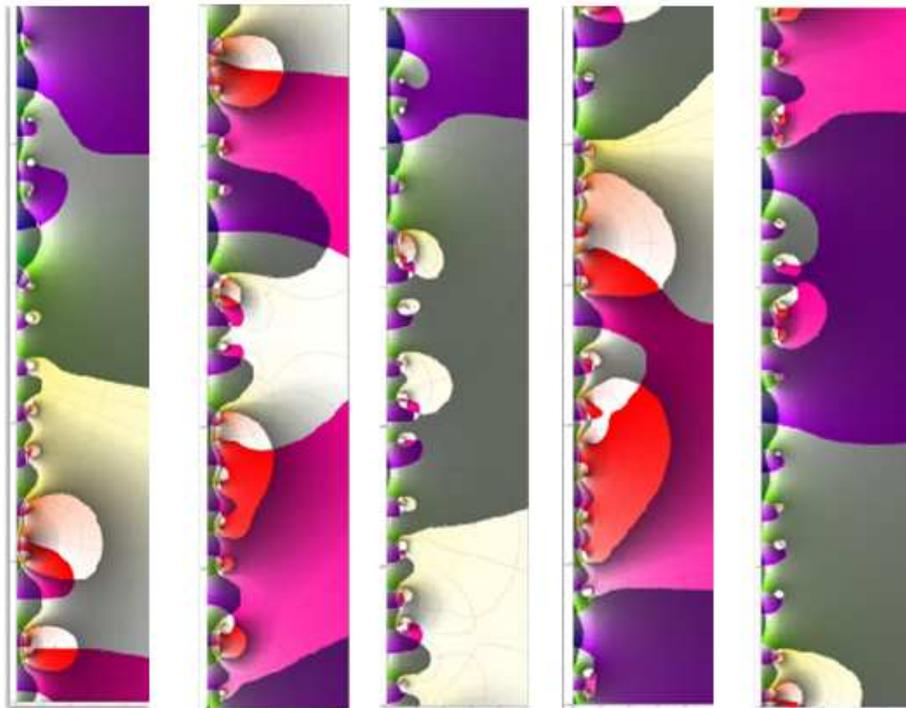

This figure exhibits two strips, one with t between one billion minus one and one billion plus 9 containing 29 non trivial zeros and the other with t between one billion plus 9 and one billion plus 17 containing 24 non trivial zeros.

Hence the number of zeros is increasing with t only on the average

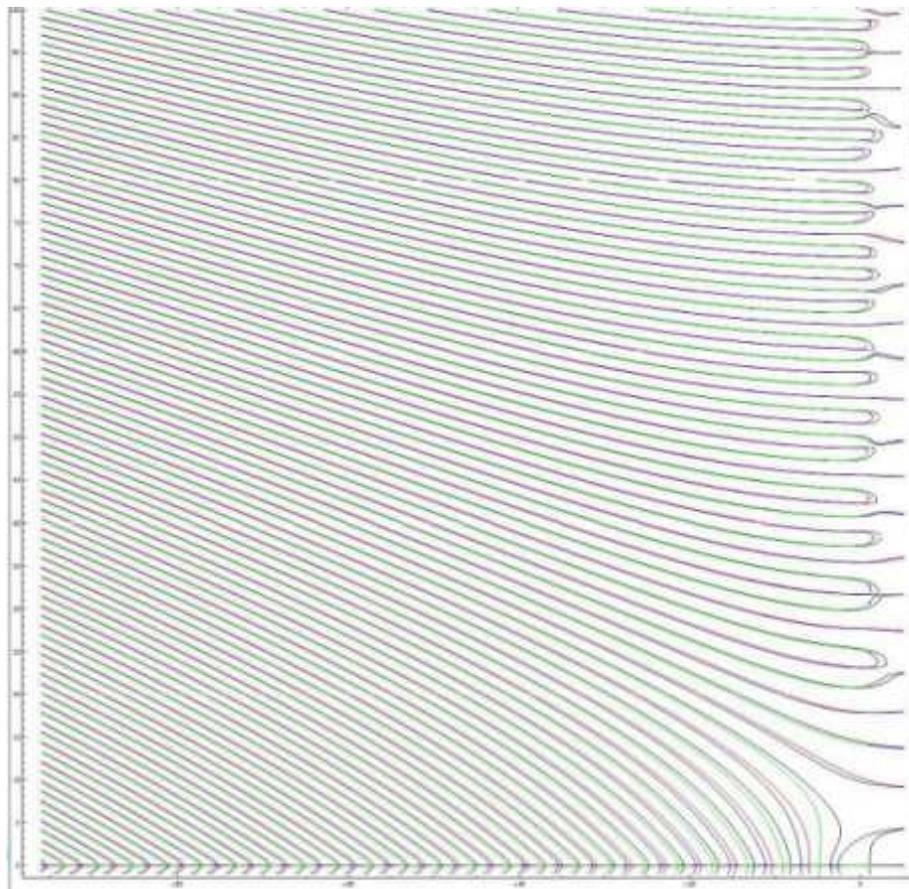

Fig. 8

This figure exhibits the intertwined curves in a box of 100X100. The curves seem to overlap two by two outside a vertical strip including the critical line.

Obviously, this is just a visual illusion, contradicting the known fact that $\Gamma_{k,j}$ and $\Upsilon_{k,j}$, or $\Gamma'_k$ and $\Upsilon'_k$ can intersect each other only in points where the tangent to $\Gamma_{k,j}$, respectively to $\Gamma'_k$ is horizontal. However, the distance on a vertical line between those curves can become smaller than the size of a printed dot, and so they appear as verlapping.

It would be an interesting project to quantize that distance by using inequalities specific to the geometric function theory. A joint paper with Grigore Salagean regarding such inequalities applied to the Riemann Zeta function has been submitted for publication.